\documentclass[10pt]{article}
\usepackage{amssymb, amsmath}
\usepackage{graphicx}
\usepackage[mathscr]{eucal}
\begin{document}
\title{The Dehn invariants of the Bricard octahedra}
\author{Victor Alexandrov}
\date{ }
\maketitle
\begin{abstract}
We prove that the Dehn invariants of any Bricard octahedron 
remain constant during the flex and that the Strong Bellows 
Conjecture holds true for the Steffen flexible polyhedron.
\end{abstract}

\renewcommand{\thefootnote}{\fnsymbol{footnote}}

\section{Introduction}
\footnotetext{The author was supported in part by the Russian State Program for Leading Scientific Schools, Grant~NSh--8526.2008.1.}

A polyhedron (more precisely, a polyhedral surface)
is said to be \textit{flexible} if its spatial shape can 
be changed
continuously due to changes of its dihedral angles only,
i.\,e., if every face remains congruent to itself during the flex.

For the first time flexible sphere-homeomorphic polyhedra
in Euclidean 3-space were constructed by R.~Connelly in
1976 \cite{Co77}. 
Since that time, many properties of flexible polyhedra were
discovered, for example: the total mean curvature \cite{Al85}
and the oriented volume \cite{Sa96}, \cite{CSW7}, \cite{Sa98}, \cite{Sc04}, 
are known to be constant during the flex.
Nevertheless, many interesting problems remain open.
One of them is the so called Strong Bellows Conjecture
which was posed in \cite{Co79} and reads as follows:
\textit{If an embedded polyhedron (i.\,e.,
self-intersection free polyhedron)
$P_1$ is obtained from an embedded 
polyhedron $P_0$ by a continuous flex
then $P_1$ and $P_0$ are scissors congruent, i.\,e.,
$P_1$ can be partitioned in a finite set of polyhedra
$\{Q_j\}$, $j=1,\dots, n$, with the following property: 
for every $j=1,\dots, n$ there is
an isometry $F_j:\mathbb R^3\to\mathbb R^3$ such that the set
$\{F_j(Q_j)\}$, is a partition of $P_0$.}

Let us recall the following well-known 

THEOREM 1. \textit{Given two embedded polyhedra $P_0$ and 
$P_1$ in $\mathbb R^3$
the following conditions are equivalent}:
\par
(1) \textit{$P_0$ and $P_1$ are scissors congruent;}
\par
(2) $\mbox{Vol\,} P_0=\mbox{Vol\,} P_1$
\textit{and} $\mbox{D}_f P_0 =\mbox{D}_f P_1$
\textit{for every $\mathbb Q$-linear function $f:\mathbb R\to\mathbb R$ such that $f(\pi)=0$.}
Here $\mbox{Vol\,} P$ stands for the volume
of $P$ and $\mbox{D}_f P=\sum |\ell| f(\alpha_\ell)$ stands for 
its Dehn invariant; $\alpha_\ell$ is the internal dihedral angle
of $P$ at the edge $\ell$; $|\ell|$ is the length of $\ell$.

The implication (1)$\Rightarrow$(2) was proved by 
M. Dehn \cite{De00}; 
the implication (2)$\Rightarrow$(1) was proved 
independently by J.P. Sydler \cite{Sy65} and  
B. Jessen \cite{Je68}.
We refer the reader to \cite{Bo78} for more detail.

Since we know that volume is constant during the 
flex \cite{Sa96}, \cite{CSW7}, \cite{Sa98}, \cite{Sc04},
Theorem~1 reduces the Strong Bellows Conjecture 
to the problem whether the Dehn invariants are constant.
Note that the latter problem makes sense also
for polyhedra with self-intersections. 
All we need is the notion of the dihedral angle 
for polyhedra with self-intersections.

DEFINITION.
A \textit{dihedral angle} at the edge $\ell$ of an oriented 
(not necessarily embedded) polyhedron $P$ 
is a multi-valued real-analytic function
$\alpha^*=\alpha_\ell^0+2\pi m$, $m\in\mathbb Z$.
Here $0<\alpha_\ell^0\leqslant 2\pi$ stands just for one of the
values (or \textit{branches}) of this multi-valued function 
and can be calculated as follows. 
Let $x$ be an internal point of the edge $\ell$,
let $Q_1$ and $Q_2$ be the two faces of $P$ that are adjacent to $\ell$,
let $\boldsymbol{n}_1$ and $\boldsymbol{n}_2$ be unit normal vectors
to the faces $Q_1$ and $Q_2$ respectively which set the orientation of $P$.
Let $\mathscr B$  be a closed ball centered at $x$ of a so small radius that~$\mathscr B$  
(a) contains no vertices of $P$, (b) has no common points with any
edge of $P$ other than $\ell$, and (c) has no common points with any
face of $P$ other than $Q_1$ or $Q_2$.
First suppose that $\boldsymbol{n}_1\not=\boldsymbol{n}_2$.
Rotate the semi-circle $\mathscr B\cap Q_1$ around its diameter 
$\mathscr B\cap\ell$
in the direction of the vector $\boldsymbol{n}_1$ until this semi-circle
coincides with the semicircle $\mathscr B\cap Q_2$  for the first time.
During the process of rotation the points of the semi-circle 
$\mathscr B\cap Q_1$ fill in a sector of $\mathscr B$. Denote this sector by $S$ and put by definition
$\alpha_\ell^0=\mbox{Vol\,} S/\mbox{Vol\,}B$.
If $\boldsymbol{n}_1=\boldsymbol{n}_2$ we put 
by definition $\alpha_\ell^0=2\pi$.

Given an edge $\ell$ of a flexible polyhedron $P(t)$
we choose an arbitrary univalent branch $\alpha_\ell (t)$ 
of the multi-valued function $\alpha_\ell^* (t)$ 
provided that $\alpha_\ell (t)$ is continuous in $t$
an use it in the calculations below.

The main result of this paper is 
that the Dehn invariants of any Bricard octahedron 
remain constant during the flex.
Using this result we prove also that the Strong Bellows 
Conjecture holds true for the Steffen flexible polyhedron.

Our description of the Bricard octahedra of types 1--3 and
of the Steffen polyhedron is very brief;
it is aimed mainly at fixing notations and recalling the properties 
needed for our study. We refer the reader for details to 
\cite{Br97}, \cite{Ku79}, \cite{Le67}, \cite{Sa92}, 
and \cite{St06}.

\section{Bricard octahedra of type 1}

Any Bricard octahedron of type 1 in $\mathbb R^3$ can be constructed in the following way.
Consider a disk-homeomorphic piece-wise linear surface $S$ in $\mathbb R^3$  composed of 
four triangles $A_1B_1C_1$, $B_1A_2C_1$, $A_2B_2C_1$, and $B_2A_1C_1$ such that
$|A_1B_1|=|A_2B_2|$ and $|B_1A_2|=|B_2A_1|$ (see Fig.~1). 
It is known  that such a spatial quadrilateral
$A_1B_1A_2B_2$ is symmetric with respect to a line~$L$ passing through the middle points 
of its diagonals $A_1A_2$ and $B_1B_2$ \cite{Ku79}.
Glue together $S$ and its symmetric image with respect to $L$ (see Fig.~2).
Denote by $C_2$ the symmetric image of~$C_1$ under the symmetry with respect to~$L$.
The resulting polyhedral surface with self-intersections is flexible (because $S$ is flexible) 
and is combinatorially equivalent to the surface of the regular octahedron. 
This is known as the Bricard octahedron of type 1.
Each of the spatial quadrilaterals $A_1B_1A_2B_2$,
$A_1C_1A_2C_2$  or $B_1C_1B_2C_2$  is called its \textit{equator}.

\begin{figure}
\includegraphics[width=0.47\textwidth]{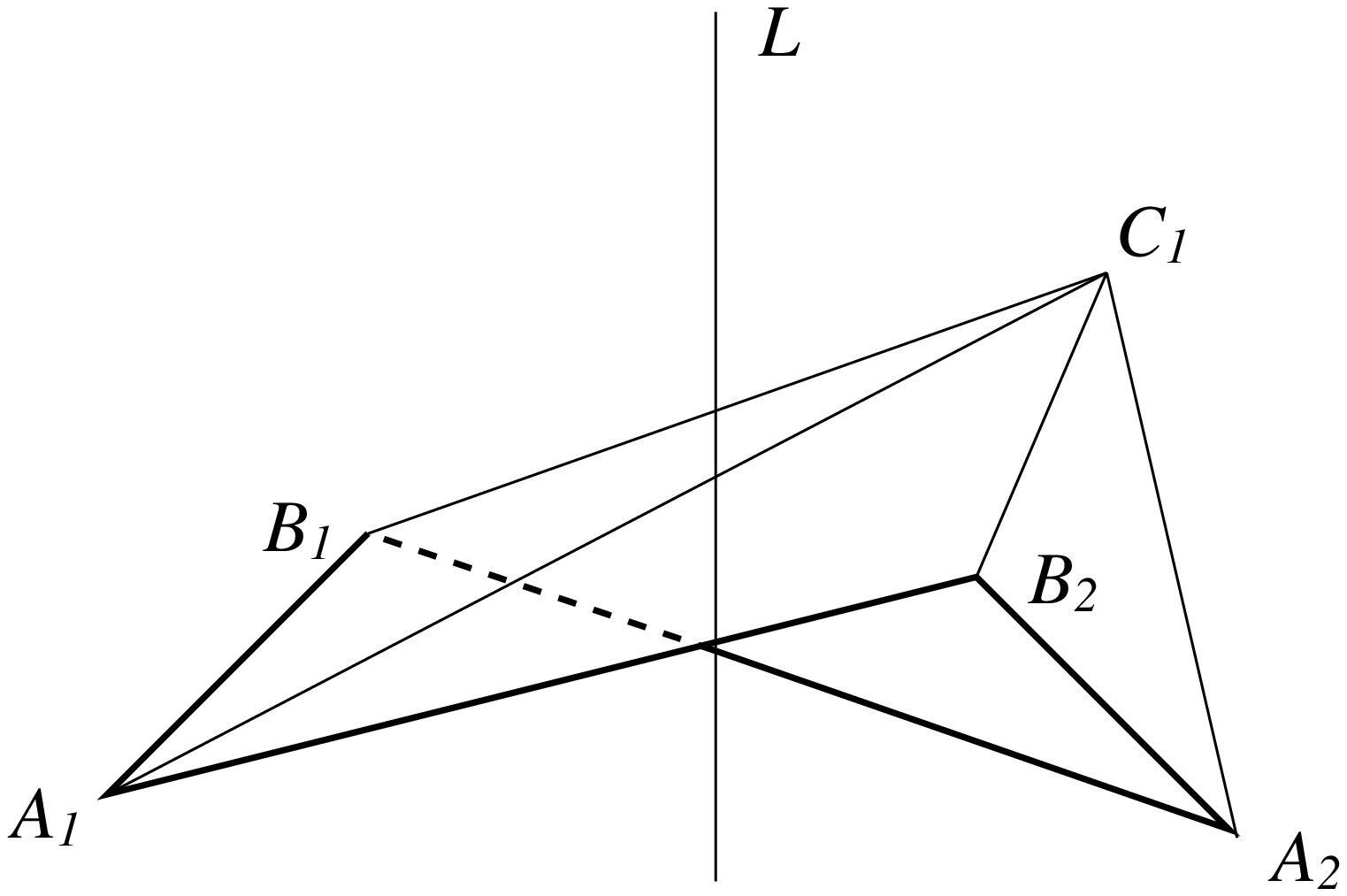}\hfill
\includegraphics[width=0.47\textwidth]{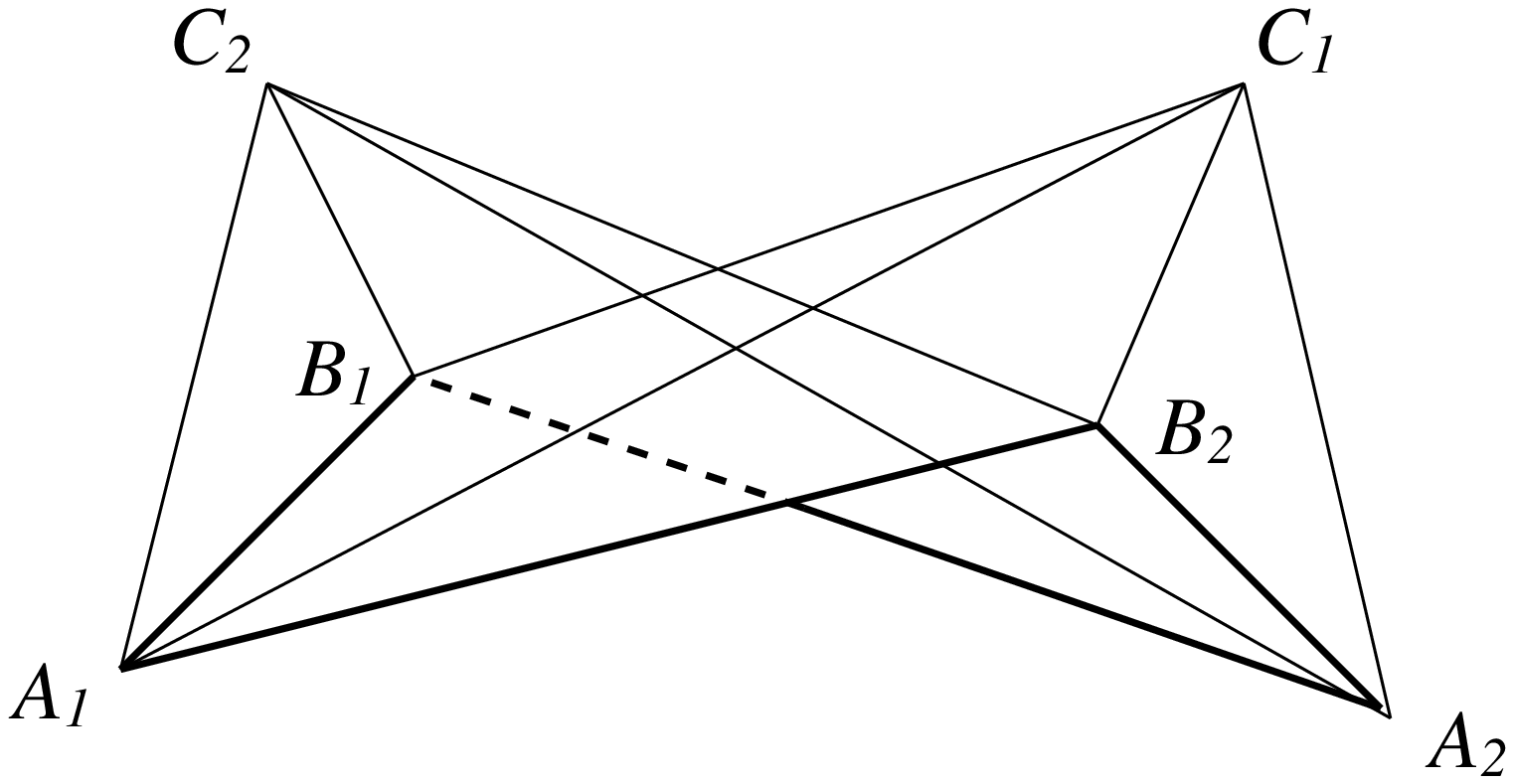}\\
\parbox[t]{0.47\textwidth}{\caption{Piece-wise linear disk $S$ \newline
                                   and line of symmetry $L$
                                   }}\hfill
\parbox[t]{0.47\textwidth}{\caption{Bricard octahedron \newline of type 1}}
\end{figure}

THEOREM 2.
\textit{For every equator $E$ of any Bricard octahedron 
of type 1 and every $\mathbb Q$-linear function 
$f:\mathbb R\to \mathbb R$ such that $f(\pi)=0$ yields}
$$
\mathscr D_f (E)\overset{\mbox{def}}{=}
\sum\limits_{\ell\in E}|\ell| f\bigl(\alpha_\ell(t)\bigr)=0.
$$

\textit{Proof}.
Split $\mathscr D_f (E)$
into 2 groups of the form
$|\ell|f\bigl(\alpha_\ell (t)\bigr)+|\ell^*|f\bigl(\alpha_{\ell^*} (t)\bigr)$,
where $\ell$ and $\ell^*$ are opposite edges of the octahedron
(i.\,e., such edges that the pair of the corresponding edges of the regular
octahedron is symmetric with respect to the center of symmetry).
Obviously, $|\ell|=|\ell^*|$. 
Moreover, $\alpha_\ell(t)+\alpha_{\ell^*}(t)=2\pi m$ for some $m\in\mathbb Z$, 
since the dihedral angles attached to $\ell$ and $\ell^*$
are symmetric to each other with respect to the line $L$.
Hence,  
$|\ell|f\bigl(\alpha_\ell (t)\bigr)+|\ell^*|f\bigl(\alpha_{\ell^*}(t)\bigr)=
|\ell|f(2\pi m)=0$ and, thus,
$\mathscr D_f (E)=0$ for all $t$. 
\hfill$\square$

THEOREM 3. 
\textit{Any Dehn invariant of any Bricard octahedron of type 1 is constant during the flex; moreover, it equals zero.}

\textit{Proof} follows immediately from Theorem 2.
\hfill$\square$

\section{Bricard octahedra of type 2}

Any Bricard octahedron of type 2 in $\mathbb R^3$ can be constructed in the following way.
Consider a disk-homeomorphic piece-wise linear surface $R$ in $\mathbb R^3$  composed of 
four triangles $A_1B_1C_1$, $B_1A_2C_1$, $A_2B_2C_1$, and $B_2A_1C_1$ such that
$|A_1B_2|=|B_2A_2|$ and $|A_1B_1|=|B_1A_2|$ (see Fig.~3). 
It is known  that such a spatial quadrilateral
$A_1B_1A_2B_2$ is symmetric with respect to a plane $P$ which dissects the
dihedral angle between the half-planes $A_1B_1B_2$ and $A_2B_1B_2$.
Glue together $R$ and its symmetric image with respect to $P$.
The resulting polyhedral surface with self-intersections is
flexible (because $R$ is flexible) and is combinatorially equivalent to 
the regular octahedron (see Fig.~4). This is known as the Bricard octahedron of type 2.
Each of the spatial quadrilaterals $A_1B_1A_2B_2$,
$A_1C_1A_2C_2$  or $B_1C_1B_2C_2$  is called its \textit{equator}.

\begin{figure}
\includegraphics[width=0.47\textwidth]{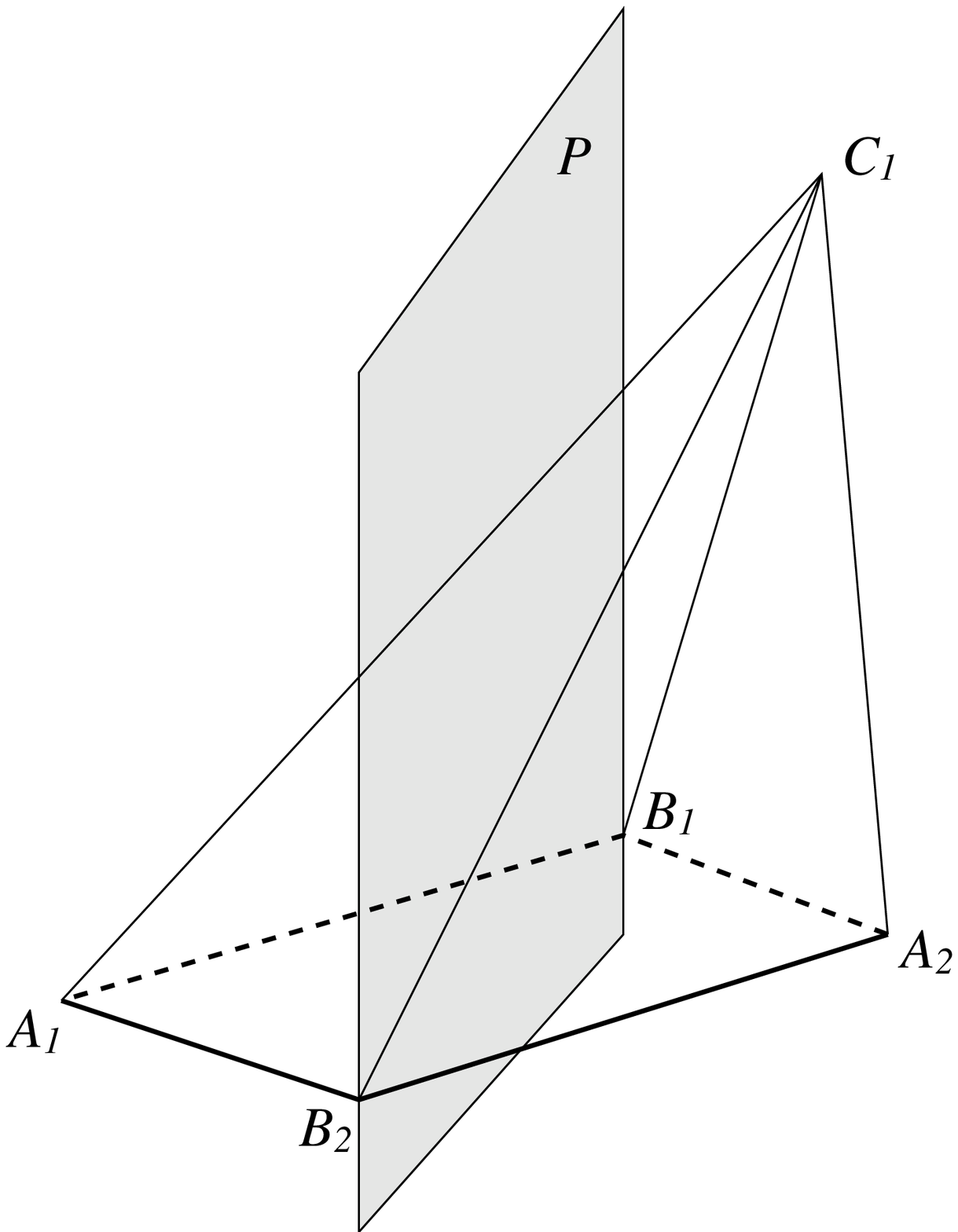}\hfill
\includegraphics[width=0.47\textwidth]{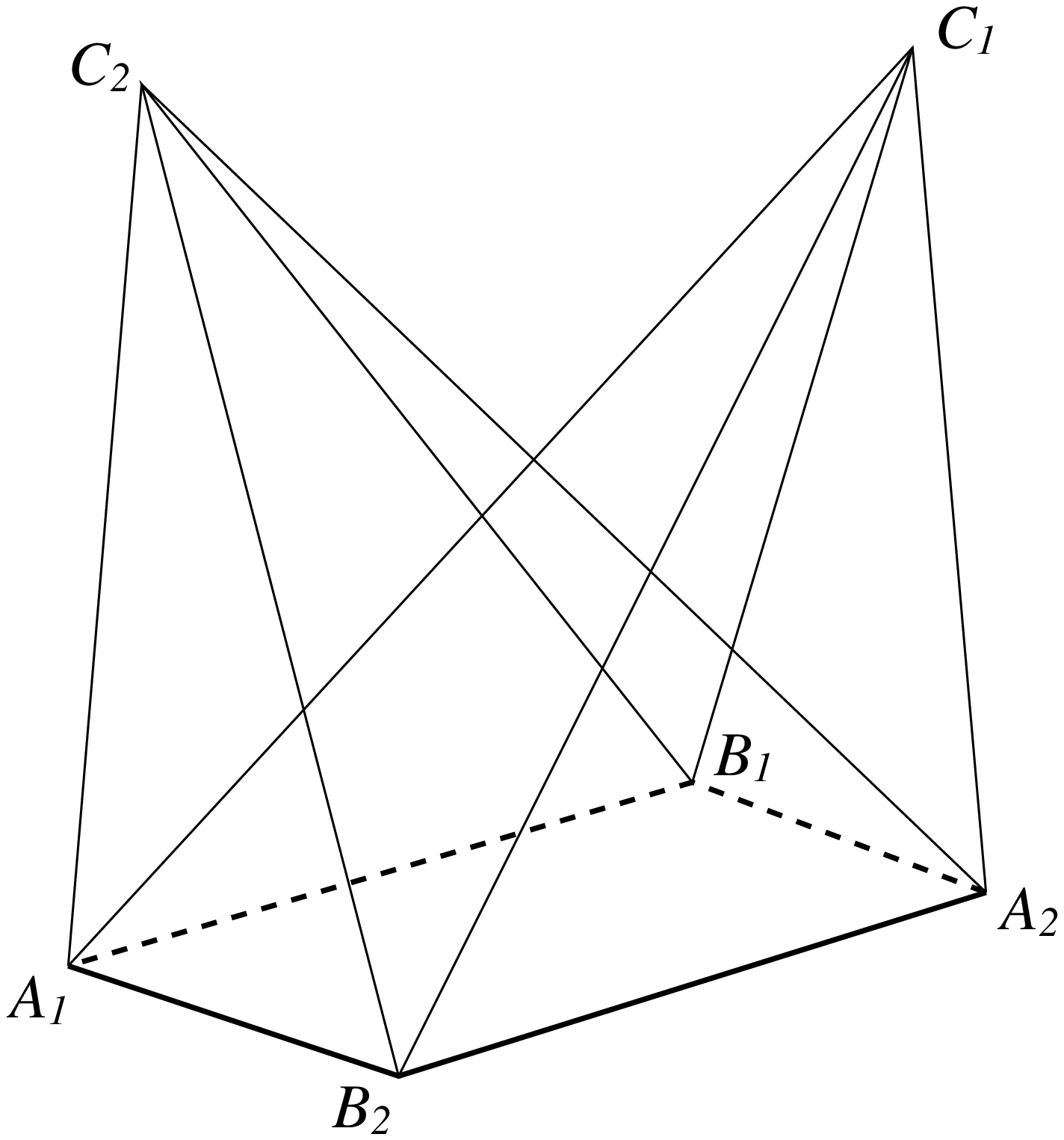}\\
\parbox[t]{0.47\textwidth}{\caption{Piece-wise linear disk $R$ \newline
                                   and plane of symmetry $P$
                                   }}\hfill
\parbox[t]{0.47\textwidth}{\caption{Bricard octahedron \newline of type 2}}
\end{figure}

THEOREM 4.
\textit{For every equator $E$ of any Bricard octahedron 
of type 2 and every $\mathbb Q$-linear function 
$f:\mathbb R\to \mathbb R$ such that $f(\pi)=0$ yields}
$\mathscr D_f (E)=0$~for~all~$t$.

\textit{Proof} is similar to the proof of Theorem 2:
symmetry of the Bricard octahedron of type 2
with respect to the plane $P$ implies that 
symmetric edges have equal lengths and 
the sum of the dihedral angles attached is a multiple of~$\pi$. 
\hfill$\square$

THEOREM 5. 
\textit{Any Dehn invariant of any Bricard octahedron 
of type 2 is constant during the flex; moreover, it equals zero.}

\textit{Proof} follows immediately from Theorem 4.
\hfill$\square$

\section{Bricard octahedra of type 3}

Any Bricard octahedron $\mathscr O$ of type 3 in $\mathbb R^3$ can be constructed in the following way.
Let $K_C$ and $K_B$ be two different circles in $\mathbb R^2$ 
with a common center.
Let $A_1B_1A_2B_2$ be a convex quadrilateral with the sides tangent to $K_C$ as it is
shown in Fig.~5
(in fact it is not prohibited that three points, say $A_1$, $B_1$, and $A_2$, lie on a straight line;
in this case $B_1$ is a tangent point).
Let $A_1C_1A_2C_2$ be a quadrilateral with self-intersections such that
every straight line containing a side of $A_1C_1A_2C_2$ is tangent to $K_B$ as it is shown in Fig.~6.
A Bricard octahedron $\mathscr O$  of type 3 in a flat position is composed of the vertices
$A_1$, $A_2$, $B_1$, $B_2$, $C_1$, and $C_2$ and of the edges of the
quadrilaterals
$A_1B_1A_2B_2$, $A_1C_1A_2C_2$, and $B_1C_1B_2C_2$ (see Fig.~7).
The faces of  $\mathscr O$  are defined as the 
triangles $\triangle A_iB_jC_k$ for all choices of 
the indices $i,j,k\in\{1,2\}$.
It is known that the quadrilateral  $B_1C_1B_2C_2$  is circumscribed about a circle
$K_A$ which shares the center with the circles $K_B$ and $K_C$ (see Fig.~7).

\begin{figure}
\includegraphics[width=0.8\textwidth]{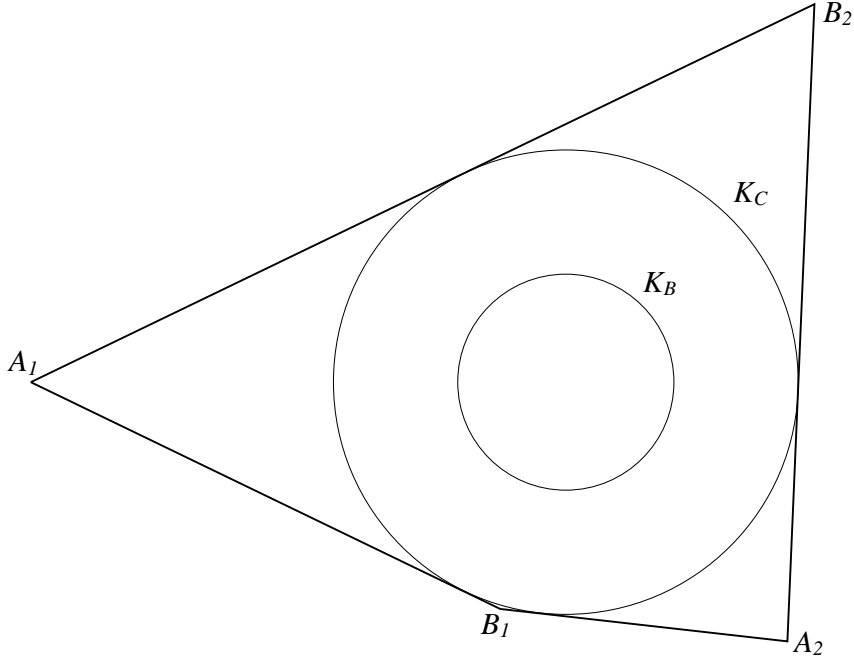}
\caption{Construction of a Bricard octahedron $\mathscr O$ of type 3. Step 1}
\end{figure}

THEOREM 6.
\textit{For every equator $E$ of any Bricard octahedron $\mathscr O$ 
of type 3 and every $\mathbb Q$-linear function 
$f:\mathbb R\to \mathbb R$ such that $f(\pi)=0$ yields}
$\mathscr D_f (E)=0$ for all $t$.

\begin{figure}
\includegraphics[width=0.8\textwidth]{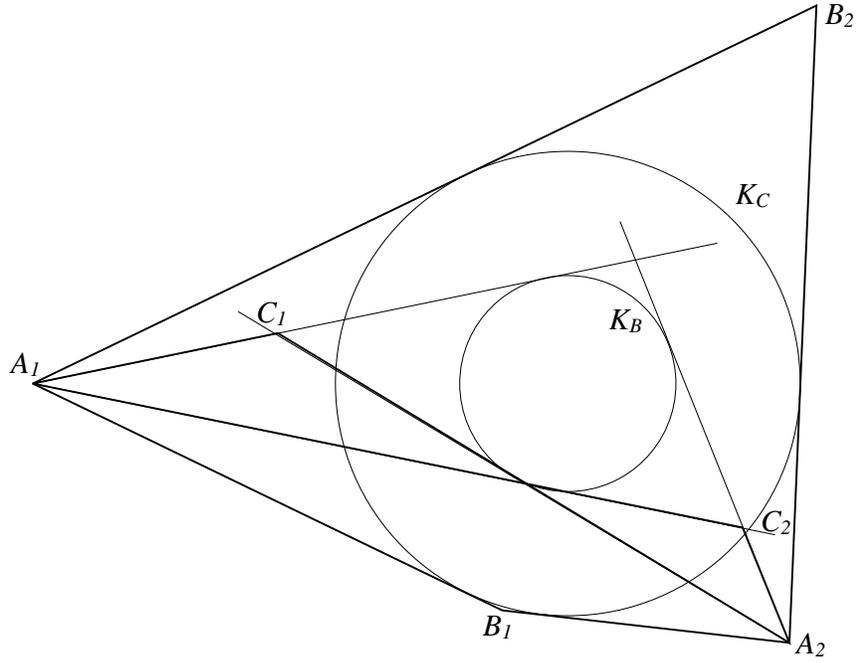}
\caption{Construction of a Bricard octahedron  $\mathscr O$ of type 3. Step 2}
\end{figure}

\textit{Proof}.
Let $E=A_1B_1A_2B_2$.  
Let $\mathscr B$ be a closed ball centered at $A_1$ 
of a so small radius~$r$
that~$\mathscr B$  contains no vertices of octahedron 
$\mathscr O$  other than~$A_1$.
The intersection of $\mathscr B$  and $\mathscr O$ is 
a quadrilateral.
Denote it by $Q(t)$ and denote its vertices as follows:
$\widetilde{B}_1=A_1B_1\cap\partial \mathscr B$,
$\widetilde{B}_2=A_1B_2\cap\partial \mathscr B$,
$\widetilde{C}_1=A_1C_1\cap\partial \mathscr B$, and
$\widetilde{C}_2=A_1C_2\cap\partial \mathscr B$.
The length of a side of $Q(t)$ equals the corresponding 
angle of a face of $\mathscr O$ multiplied by the 
radius $r$ of $\mathscr B$ , e.\,g.,
the length of the side $\widetilde{B}_2\widetilde{C}_1$
is equal to $\angle B_2A_1C_1$ multiplied by $r$.
Thus, the length of any side of $Q(t)$ remains 
constant during the flex.
On the other hand, it follows from Fig.~7 that 
$\angle B_2A_1C_1=\angle B_1A_1C_2$ and
$\angle B_2A_1C_2=\angle B_1A_1C_1$.
Hence, the opposite sides of $Q(t)$ have equal lengths.
\begin{figure}
\includegraphics[width=0.8\textwidth]{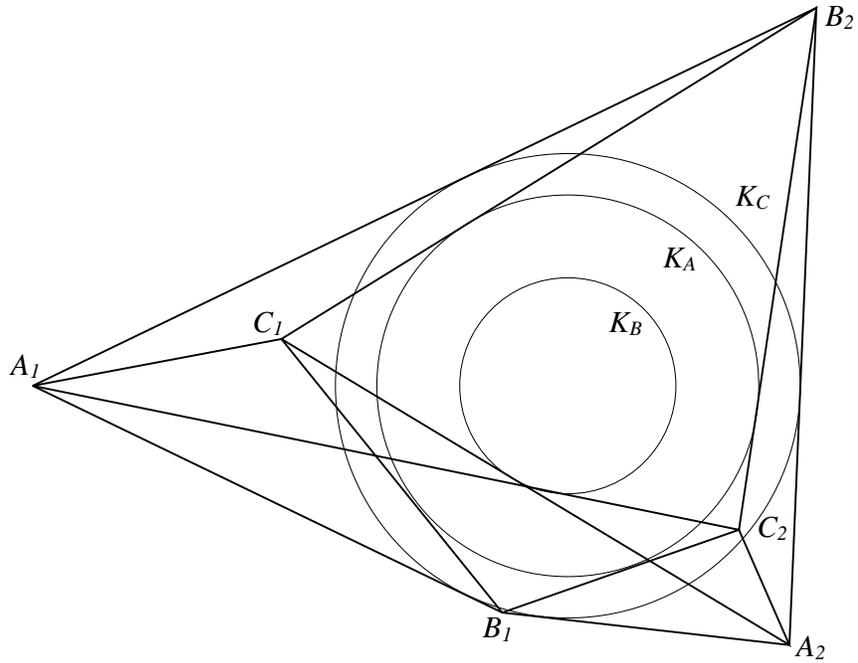}
\caption{Construction of a Bricard octahedron  $\mathscr O$ of type 3. Step 3}
\end{figure}

If $\mathscr O$  is in the flat position shown in Fig.~7 then
all the vertices of the spherical quadrilateral $Q(t)$ are 
located on a single great circle (see Fig. 8a). 
If $\mathscr O$ is in a non-flat position close to a flat position
shown in Fig.~7 then, generally speaking, there are two 
possibilities for $Q(t)$: either $Q(t)$ is convex (see Fig.~8b)
or $Q(t)$ has self-intersections (see Fig.~8c).
However, $Q(t)$ can not be convex because in this case
$\mathscr O$ is embedded (at least when $\mathscr O$ 
is close enough
to the flat position shown in Fig.~7) while it is known that
no embedded  octahedron is flexible  \cite{Ma08}.
Taking into account that opposite sides of $Q(t)$ have equal
lengths we conclude that the sum of any two opposite angles 
of $Q(t)$ equals $2\pi m$ for some $m\in\mathbb Z$.
Hence, the  sum of any two opposite dihedral angles 
of the solid angle of $\mathscr O$ with the vertex
$A_1$ equals $2\pi m$ for some $m\in\mathbb Z$.
According to our agreement made in the Introduction,
this means, for example, that if $\alpha_{A_1B_1} (t)$
and $\alpha_{A_1B_2} (t)$ are any branches of the dihedral
angles of $\mathscr O$ at the edges $A_1B_1$ and $A_1B_2$,
respectively, then there is $n\in\mathbb Z$ such that
for all $t$ yields
$$
\alpha_{A_1B_1} (t) +\alpha_{A_1B_2} (t)=2\pi n.  \eqno(1)
$$

Strictly speaking, the above arguments prove (1) only for 
the positions of~$\mathscr O$ close enough to the flat position
shown in Fig.~7. But, in fact, (1) holds true
for all positions of $\mathscr O$ obtained from the flat
position shown in Fig.~7 by a continuous flex.
The reason is that we may assume that the coordinates 
of the vertices of $\mathscr O$ are analytic functions 
of the flexing parameter $t$, see \cite{Gl75} for more detail.
Then any branch of the dihedral angle is an analytic function 
of $t$.
As soon as we know that (1) holds true
for all $t$ corresponding to any position of~$\mathscr O$
close enough to the flat position shown in Fig.~7,
we conclude that
(1) holds true for all $t$ corresponding to any position 
of $\mathscr O$ obtained from the flat
position shown in Fig.~7 by a continuous flex.
But the configuration space (i.\,e., the
space of all positions)  of any Bricard octahedron 
of type 3 is known to be homeomorphic to a circle.
Thus the relation between dihedral angles established
for $\mathscr O$ close enough to the flat position
shown in Fig.~7 holds true for all positions of $\mathscr O$.

Similar arguments show that, for 
any solid angle of $\mathscr O$ with the vertex
$B_1$, $A_2$, or $B_2$,
the  sum of any two opposite dihedral angles 
equals $2\pi m$ for some $m\in\mathbb Z$.

By definition put $A_iB_j\cap K_C=c_{ij}$; $i,j=1,2$.
In other words, denote by $c_{ij}$ the point where the segment
$A_iB_j$ touches the circle $K_C$.
Then
\begin{align*}
\mathscr D_f (E)=&\bigl(|A_1c_{12}|f(\alpha_{A_1B_2})+
|A_1c_{11}|f(\alpha_{A_1B_1})\bigr)\\
&+\bigl(|B_1c_{11}|f(\alpha_{A_1B_1})+
|B_1c_{21}|f(\alpha_{A_2B_1})\bigr)\\
&+\bigl(|A_2c_{21}|f(\alpha_{A_2B_1})+
|A_2c_{22}|f(\alpha_{A_2B_2})\bigr)\\
&+\bigl(|B_2c_{22}|f(\alpha_{A_2B_2})+
|B_2c_{12}|f(\alpha_{A_1B_2})\bigr)\\
=&|A_1c_{12}|f(\alpha_{A_1B_2}+\alpha_{A_1B_1})
 +|B_1c_{11}|f(\alpha_{A_1B_1}+\alpha_{A_2B_1})\\
&+|A_2c_{21}|f(\alpha_{A_2B_1}+\alpha_{A_2B_2})
 +|B_2c_{22}|f(\alpha_{A_2B_2}+\alpha_{A_1B_2})=0.
\end{align*}
Here we use the fact that
the segments of the two lines passing through 
a point and tangent to a circle have
equal lengths (e.\,g., $|A_1c_{12}|=|A_1c_{11}|$)
and that the sum of the opposite angles
of $Q(t)$ equals $2\pi m$ for some $m\in\mathbb Z$
(e.\,g., (2)).

\begin{figure}
\includegraphics[width=0.37\textwidth]{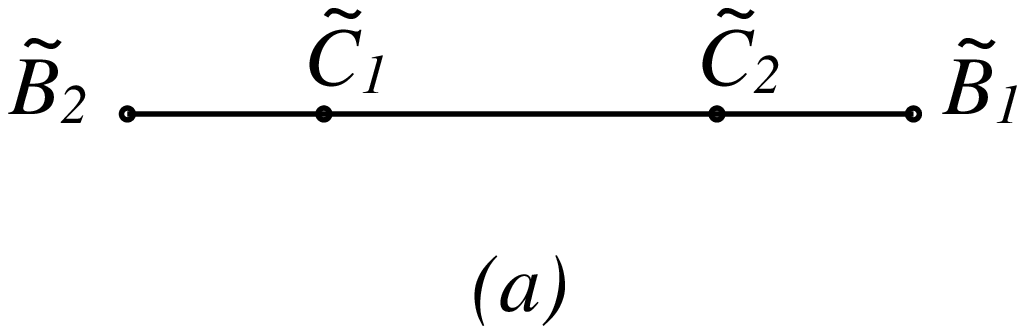}
\quad\includegraphics[width=0.25\textwidth]{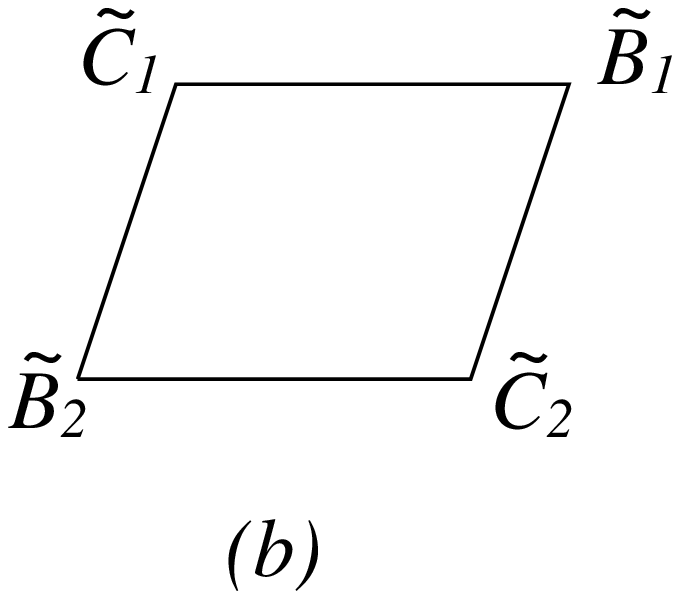}\hfill
\includegraphics[width=0.3\textwidth]{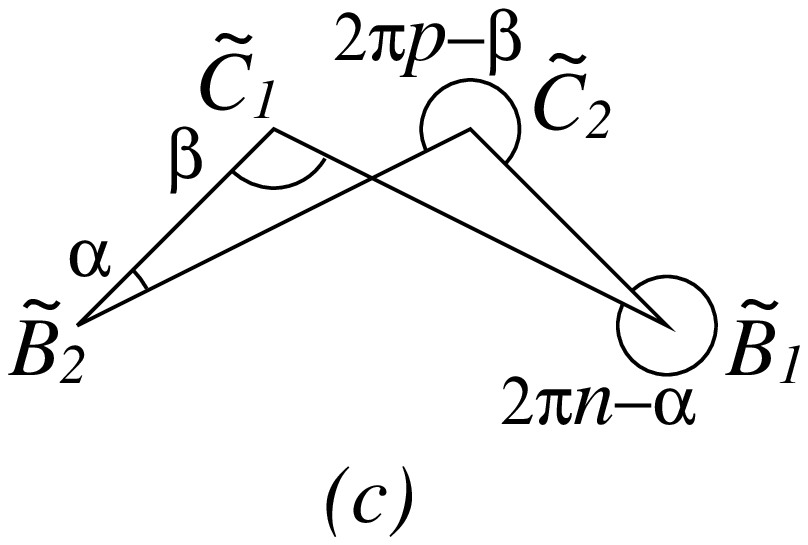}\quad\\
\caption{Spherical quadrilateral $Q(t)$: (a) in a flat position;
(b) as a convex polygon (impossible); (c) as a polygon a with 
self-intersection}
\end{figure}

Exactly the same arguments, we have used above to prove 
the formula $\mathscr D_f (E)=0$ for the equator $E=A_1B_1A_2B_2$, 
can be used to prove it for the equator $E=B_1C_1B_2C_2$.
This part of the proof is left to the reader. 
On the contrary, in order to prove the formula 
$\mathscr D_f (E)=0$ for the equator $E=A_1C_2A_2C_1$, 
we have to use additional arguments that are given below.

Let $E=A_1C_2A_2C_1$.  
We need some relations for the dihedral angles
at the edges of the equator~$E$. 
We already know that if $\alpha=\alpha(t)$ is 
one of the univalent branches of the angle $\alpha_{A_1B_2}(t)$
(i.\,e., $\alpha_{A_1B_2}(t)=\alpha(t)$) then
$\alpha_{A_2B_2}(t)=2\pi k-\alpha$,
$\alpha_{A_2B_1}(t)=\alpha +2\pi m$, and
$\alpha_{B_1A_1}(t)=2\pi n -\alpha$
for some $k,m,n\in\mathbb Z$.
Similarly, we know that 
if $\beta=\alpha_{A_1C_1}(t)$ and
$\gamma=\alpha_{C_2A_2}(t)$ then
$\alpha_{A_1C_2}(t)+\beta=2\pi p$ and 
$\alpha_{A_2C_1}(t)+\gamma=2\pi q$
for some $p,q\in\mathbb Z$.
Let's prove that, for some $s\in\mathbb Z$, yields
$$\beta(t)+\gamma(t)=2\pi s.\eqno(2)$$

Recall one of the four Napier's analogies 
(known also as the Napier's rules)
from the spherical trigonometry \cite{Sm60}:
\textit{Let a spherical triangle have sides $a$, $b$, and $c$ 
with $A$, $B$, and $C$ the corresponding opposite angles. Then}
$$
\frac{\sin\dfrac{a-b}{2}}{\sin\dfrac{a+b}{2}}=
\frac{\tan\dfrac{A-B}{2}}{\cot\dfrac{C}{2}}.\eqno(3)
$$

Applying  (3) to the spherical triangle
$\widetilde{C}_1\widetilde{B}_2\widetilde{C}_2$
shown in Fig.~8c we get
$$
\frac{\sin r\dfrac{\angle B_2A_1C_2-\angle B_2A_1C_1}{2}}{\sin r\dfrac{\angle B_2A_1C_2+\angle B_2A_1C_1}{2}}=
\frac{\tan\dfrac{\beta}{2}}{\cot\dfrac{\alpha}{2}},\eqno(4)
$$
where $r$ stands for the radius of the ball $\mathscr B$
and $\angle XYZ$ stands for the value of the `plane' angle 
with the vertex $Y$ of the triangular
face $XYZ$ of the Bricard octahedron $\mathscr O$
shown in Fig.~7.

Similarly, apply (3) to a spherical triangle obtained by
intersecting the solid angle of the octahedron $\mathscr O$
with the vertex $A_2$ and the ball $\mathscr B^*$ of the radius
$r$ centered at $A_2$ (we assume that $\mathscr B^*$ contains
no vertices of $\mathscr O$ other than $A_2$ and the
balls $\mathscr B$, $\mathscr B^*$ have equal radii).
More precisely, put
$\widetilde{B}_1^*=A_2B_1\cap\partial \mathscr B^*$,
$\widetilde{C}_1^*=A_2C_1\cap\partial \mathscr B^*$, and
$\widetilde{C}_2^*=A_2C_2\cap\partial \mathscr B^*$.
Applying (3) to the spherical triangle 
$\widetilde{B}_1^*\widetilde{C}_1^*\widetilde{C}_2^*$
(see Fig.~9) we get
$$
\frac{\sin r\dfrac{\angle B_1A_2C_2-\angle B_1A_2C_1}{2}}{\sin r\dfrac{\angle B_1A_2C_2+\angle B_1A_2C_1}{2}}=
\frac{\tan\biggl(\pi q-\dfrac{\gamma}{2}\biggr)}{\cot\biggl(\dfrac{\alpha}{2}+\pi m\biggr)}.\eqno(5)
$$
\begin{figure}
\begin{center}
\includegraphics[width=0.37\textwidth]{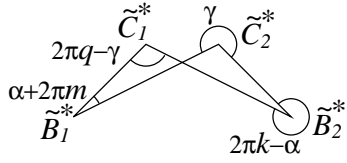}
\end{center}
\caption{Spherical quadrilateral 
$\widetilde{B}_1^*\widetilde{C}_1^*\widetilde{B}_2^*\widetilde{C}_2^*$}
\end{figure}

Using Fig.~6 we easily conclude that 
$\angle B_2A_1C_2=\angle B_1A_2C_2$ and
$\angle B_2A_1C_1=\angle B_1A_2C_1$.
Hence, the left-hand sides of (4) and (5)
are equal to each other implying that 
$\beta+\gamma=2\pi s$
for some $s\in\mathbb Z$.
Thus, the formula (2) is proved.

By definition put $A_iC_j\cap K_B=b_{ij}$; $i,j=1,2$.
In other words, denote by $b_{ij}$ the point where the 
line passing through the segment
$A_iC_j$ touches the circle $K_B$.
Then $|A_1C_1|=|A_1b_{11}|-|C_1b_{11}|$,
$|C_1A_2|=|C_1b_{21}|+|b_{21}A_2|$,
$|A_2C_2|=|A_2b_{22}|-|C_2b_{22}|$, and
$|C_2A_1|=|C_2b_{12}|+|b_{12}A_1|$.
Besides, $|A_1b_{11}|=|A_1b_{12}|$,
$|C_1b_{11}|=|C_1b_{21}|$,
$|C_2b_{12}|=|C_2b_{22}|$, and
$|A_2b_{21}|=|A_2b_{22}|$.

Now we can make the main computation for the equator
$E=A_1C_2A_2C_1$: 
\begin{align*}
\mathscr D_f (E)=& {}\bigl(|A_1b_{11}|f(\alpha_{A_1C_1})+ 
|A_1b_{12}|f(\alpha_{A_1C_2})\bigr)\\  
&-\bigl(|C_1b_{11}|f(\alpha_{A_1C_1})- 
|C_1b_{21}|f(\alpha_{A_2C_1})\bigr)\\  
&+\bigl(|A_2b_{21}|f(\alpha_{A_2C_1})+ 
|A_2b_{22}|f(\alpha_{A_2C_2})\bigr)\\  
&-\bigl(|C_2b_{22}|f(\alpha_{A_2C_2})- 
|C_2b_{12}|f(\alpha_{A_1C_2})\bigr)\\
=& {}|A_1b_{11}|f(\alpha_{A_1C_1}+\alpha_{A_1C_2})    
 +|C_1b_{21}|f(\alpha_{A_2C_1}-\alpha_{A_1C_1})\\  
&+|A_2b_{22}|f(\alpha_{A_2C_1}+\alpha_{A_2C_2})    
 +|C_2b_{12}|f(\alpha_{A_1C_2}-\alpha_{A_2C_2})\\
=& {}|A_1b_{11}|f(2\pi p)    
 +|C_1b_{21}|f(2\pi q-2\pi s)\\  
&+|A_2b_{22}|f(2\pi q)    
 +|C_2b_{12}|f(2\pi p-2\pi s)=0. 
\qquad\qquad\qquad\square
\end{align*}
 
THEOREM 7. 
\textit{Any Dehn invariant of any Bricard octahedron of type 3 is constant during the flex; moreover, it equals zero.}

\textit{Proof} follows immediately from Theorem 6. \hfill $\square$

REMARK.
Bricard octahedra of type 3 are unexpectedly symmetric: 
we have seen above that there are linear relations between
edge lengths, plane and dihedral angles, and
trigonometric relations between dihedral angles.
Let us now mention one more relation that did not appear
above: $\angle A_1C_1B_2+\angle A_2C_1B_1=
\angle A_1C_1B_1+\angle A_2C_1B_2=\pi$ (see Fig.~7).
A proof is left to the reader.

\section{Steffen polyhedron}

There are several examples of embedded flexible polyhedra
in $\mathbb R^3$ proposed by R. Connelly  (see \cite{Co77} or \cite{Sa02}),
by P. Deligne and N. Kuiper (see \cite{Sa02}), 
and by K. Steffen (see \cite{Be87} or \cite{Sa02}).
Each of them belongs to the class $\mathscr{F}_n$ for some $n\geqslant 0$,
where $\mathscr{F}_n$ is defined as follows:

(i) $\mathscr{F}_0$ consists of all convex polyhedra in $\mathbb R^3$  and all
Bricard octahedra;

(ii) a polyhedron $P$ belongs to $\mathscr{F}_n$, $n\geqslant 1$,
if and only if one of the following holds true: 
(ii${}_1$) $P$ is obtained from  $P_1, P_2\in\mathscr{F}_k$, $0\leqslant k\leqslant n-1$, 
by gluing them together along congruent faces
$Q_1\subset P_1$ and $Q_2\subset P_2$;
(ii${}_2$) $P$ is obtained from $P_1\in\mathscr{F}_{n-1}$ by gluing together two its faces
$Q_1,Q_2\subset P_1$ provided that they coincide in $\mathbb R^3$;
(ii${}_3$) $P$ is obtained from $P_1\in\mathscr{F}_{n-1}$ by a subdivision of its faces.

For example, the Steffen polyhedron (which has only 9 vertices 
and 14 faces and, hypothetically, is an embedded flexible polyhedron
with the least possible number of vertices)
can be constructed from a tetrahedron $\mathscr T$ and two copies $\mathscr O$ and $\mathscr O^\dag$
of the Bricard octahedron of type 1 in the following way (see Fig.~10).
 
\begin{figure}
\begin{center}
\includegraphics[width=0.56\textwidth]{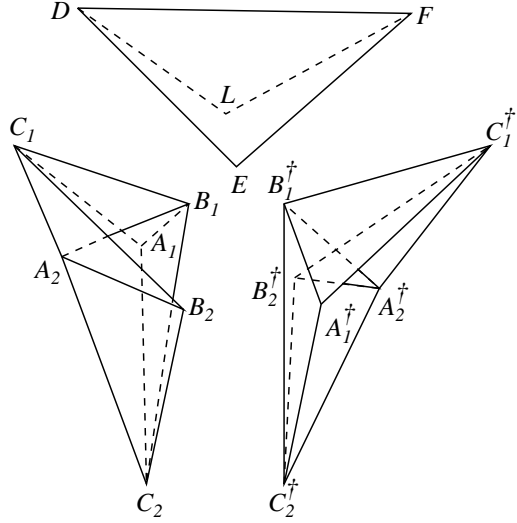}
\end{center}
\caption{Constructing the Steffen polyhedron}
\end{figure}

The tetrahedron $\mathscr T=DEFL$ has the following edge lengths:
$|DE|=|EF|=|FL|=|LD|=12$ and $|DF|=17$. It dos not change its spatial shape
during the flex of the Steffen polyhedron. The edge $EL$ is not shown in
Fig.~10 because it does not appear in the Steffen polyhedron.

The Bricard octahedron  of type 1 $\mathscr O=A_1A_2B_1B_2C_1C_2$
has the following edge lengths:
$|A_1C_1|=|C_2B_2|=|A_2C_2|=|C_2B_1|=12$,
$|B_1C_1|=|C_1A_2|=|A_1C_2|=|C_2B_2|=10$,
$|A_1B_1|=|A_2B_2|=5$, and
$|A_1B_2|=|A_2B_1|=11$.
The edge $A_1B_2$ is not shown in
Fig.~10 because it does not appear in the Steffen polyhedron.

The Bricard octahedron  of type 1 $\mathscr O^\dag =A_1^\dag A_2^\dag B_1^\dag B_2^\dag C_1^\dag C_2^\dag$
is obtained from $\mathscr O$ by an orientation-preserving isometry of $\mathbb R^3$.
The edge $A_1^\dag B_2^\dag$ is not shown in
Fig.~10 because it does not appear in the Steffen polyhedron.

Glue $\mathscr T$ and $\mathscr O$ along the triangles $\triangle DEL$ and $\triangle C_1B_2A_1$
(more precisely, identify the points $D$ and $C_1$, $E$ and $B_2$, $L$ and $A_1$).
The resulting polyhedron $\mathscr S_1$ belongs to the class $\mathscr F_1$, is flexible
but is not embedded. 

Glue $\mathscr P_1$ and $\mathscr O^\dag$ along the triangles $\triangle EFL$ and $\triangle A_1^\dag C_1^\dag B_2^\dag$
(more precisely, identify the points $E$ and $A_1^\dag$, $F$ and $C_1^\dag$, $L$ and $B_2^\dag$).
The resulting polyhedron $\mathscr S_2$ belongs to the class $\mathscr F_2$, is flexible
but is not embedded. 

Note that, during the flex of $\mathscr S_2$, the both vertices $C_2$ and $C_2^\dag$ move along the circle 
that lies in the plane perpendicular to the segment $EL$ and is centered at the middle point of $EL$.
Hence, for every position of $C_2$ (originated from the flex of $\mathscr O$) we can bend $\mathscr O^\dag$
in such a way that $C_2^\dag$  coincide with $C_2$. It means that even when we glued the triangles
$\triangle LEC_2$ and $\triangle LEC_2^\dag$ the resulting polyhedron is flexible.
It belongs to the class $\mathscr F_3$ and is known as the Steffen flexible polyhedron, see \cite{Be87} or \cite{Sa92}.

THEOREM 8. \textit{For every $n\geqslant 0$ any flexible embedded polyhedron $P\in\mathscr{F}_n$
satisfies the Strong Bellows Conjecture, i.\,e., every embedded polyhedron $P'$  obtained
from $P$ by a continuous flex is scissors congruent to $P$.}

\textit{Proof.} According to Theorem 1 it suffice to prove that
$\mbox{Vol\,} P'=\mbox{Vol\,} P$ and 
$$\mbox{D}_f P' =\mbox{D}_f P\eqno(6)$$
for every $\mathbb Q$-linear function $f:\mathbb R\to\mathbb R$ such that $f(\pi)=0$. 
The equality $\mbox{Vol\,} P'=\mbox{Vol\,} P$ follows directly from the fact that
the oriented volume is constant during the flex \cite{Sa96}, \cite{CSW7}, \cite{Sa98}, \cite{Sc04}.
The equality (6) can be proven by induction on~$n$. 
In fact, Theorems 3, 5, and 7 yields (6) for $P\in\mathscr{F}_0$.
Now suppose that (6) holds true for all polyhedra of the classes $\mathscr{F}_k$, $0\leqslant k\leqslant n-1$.
If $P$ is constructed according to (ii${}_1$) from $P_1, P_2\in\mathscr{F}_k$, $0\leqslant k\leqslant n-1$, then
$\mbox{D}_f P =\pm\mbox{D}_f P_1  \pm\mbox{D}_f P_2 $ (depending on the orientation of 
$P$, $P_1$, and $P_2$) and, thus, is constant during the flex.
Similarly, if $P$ is constructed according to (ii${}_2$) or (ii${}_3$) from 
$P_1\in\mathscr{F}_{n-1}$ then $\mbox{D}_f P =\pm\mbox{D}_f P_1$ and, thus, 
is again constant during the flex.
\hfill $\square$

COROLLARY. \textit{If $P$ is a Steffen flexible polyhedron
and $P'$ is obtained from $P$ by a continuous flex then $P$ and $P'$ are scissors 
congruent.}\hfill $\square$

\bigskip

\noindent\textit{Victor Alexandrov}

\noindent\textit{Sobolev Institute of Mathematics}

\noindent\textit{Koptyug ave., 4}

\noindent\textit{Novosibirsk, 630090, Russia}

\noindent and

\noindent\textit{Department of Physics}

\noindent\textit{Novosibirsk State University}

\noindent\textit{Pirogov str., 2}

\noindent\textit{Novosibirsk, 630090, Russia}

\noindent\textit{e-mail: alex@math.nsc.ru}

\bigskip

\noindent{Received January 20, 2009}

\end{document}